%% file: Hilbertandbidisk.tex
\documentclass{birkmult}

\usepackage{latexsym}
\usepackage{amssymb}

\let\cal\mathcal

\input{macros_withbirk}

\renewcommand{\theequation}{\thesection.\arabic{equation}}

\newcommand{\G}{\Gamma}
\renewcommand{\z}{\zeta}

\newcommand{\hido}{H_1^{\infty}(\D)}
\newcommand{\hidd}{H_1^{\infty}(\D^d)}
\newcommand{\hib}{\hb}
\newcommand{\hibo}{H^{\i}_1(\D^2)}
\renewcommand{\E}{{\mathbb E}}
\renewcommand{\B}{{\mathbb B}}
\renewcommand{\O}{\Omega}

\begin{document}
\setlength{\baselineskip}{21pt}
\title{What Hilbert spaces can 
tell us about bounded functions in the bidisk}
\author{Jim Agler}
\address{ U.C. San Diego,
La Jolla, California 92093}
\author{John E. M\raise.5ex\hbox{c}Carthy}
\address{ Washington University,
St. Louis, Missouri 63130}
\thanks{Jim Agler was partially supported by National Science Foundation Grant
DMS 0400826, and John  M\raise.5ex\hbox{c}Carthy was
partially supported by National Science Foundation Grant
DMS 0501079.}
\dedicatory{Dedicated to the memory of Paul Halmos}

\bibliographystyle{plain}

\begin{abstract}
We discuss various theorems about bounded analytic functions on the bidisk that were proved using
operator theory. 
\end{abstract}

\maketitle

\baselineskip = 18pt

\section{Introduction}\label{seca}
Much of modern single operator theory (as opposed to the study of operator algebras)
rests on a foundation of complex analysis. Every cyclic operator can be represented as multiplication
by the independent variable on the completion of the polynomials with respect to some norm.
The nicest way to have a norm is in $L^2(\mu)$, and then one is led to the study
of subnormal operators, introduced by P.~Halmos in \cite{hal}. The study of cyclic subnormal operators
becomes the study of the spaces $P^2(\mu)$, the closure of the polynomials in $L^2(\mu)$, and 
the theory of these spaces relies on a blend of complex analysis and functional analysis; see
J.~Conway's book \cite{Co91} for an exposition. Alternatively, one can start with a Hilbert space
that is amenable to complex analysis, such as the Hardy space $H^2$, and study classes of operators
on that space that have a good function theoretic representation, such as Toeplitz, Hankel or composition operators.
All of these classes of operators have a rich theory, 
which depends heavily on function theory --- for expositions, see \eg
\cite{nik02, pel02} and \cite{cm95}.

The traffic, of course, goes both ways. There are many questions in function theory that have either
been answered or illuminated by an operator theory approach. The purpose of this article is to describe
how operator theory has fared when studying $\hb$, the algebra of bounded 
analytic functions on the bidisk $\Dt$.
We focus on function theory results that were proved, originally at least,
using operator theory.

For the topics in Sections~\ref{secb} to \ref{secff}, 
we shall first describe the situation on the disk $\D$, and then move on to the bidisk.
The topics in Sections~\ref{sech} to \ref{secj} do not really have analogues in one dimension.
For simplicity, we shall stick to scalar-valued function theory,
though many of the results have natural matrix-valued analogues.

We shall use the notation that points in the bidisk are called $\l$ or $\z$, and
their coordinates will be given by superscripts: $ \l = (\l^1,\l^2 )$.
We shall use $z$ and $w$ to denote the coordinate functions on $\D^2$.
The closed unit ball of $\hid$ will be written $\hido$ and the closed unit ball of
$\hib$ as $\hibo$.

\section{Realization Formula}
The realization formula is a way of associating isometries (or contractions) with functions
in the ball of $\hid$. In one dimension, it looks like the following; see \eg \cite{bgr90} or \cite{ampi} 
for a proof.
\label{secb}
\bt
\label{thmec1}
The function $\phi$ is in the closed unit ball of $\hid$ if and
only if there is a Hilbert space $\h$ and an isometry 
$V : \C \oplus \h \to \C \oplus \h$, such that, writing
$V$ as
\be
\label{eqecb1}
V \= 
\bordermatrix{&\C &\h \cr
\C &A & B \cr
\h &C  & D} ,
\ee
one has
\be
\label{eqecb2}
\phi(\l) \= A + \l B(I-\l D)^{-1}C .
\ee
\et
This formula was generalized to the bidisk in \cite{ag90}.
It becomes
\bt
\label{thmka2}
The function $\phi$ is in the closed unit ball of $H^\i(\D^2)$
if and only
if there are auxiliary Hilbert spaces $\h_1$ and $\h_2$ and an
isometry $$
V\, :\,\C \oplus \h_1 \oplus \h_2
\ \to \
\C \oplus \h_1 \oplus \h_2
$$ such that, if $\h := \h_1 \oplus
\h_2$, $V$ is written as
\be
\label{eqka22}
V \= \bordermatrix{&\C &\h\cr
\C &A & B\cr
\h&C &D} ,
\ee
and $\Ee = \l^1 I_{\h_1} \oplus \l^2 I_{\h_2} $,
then
\be
\label{eqka3}
\phi(\lambda) \= A + B \Ee (I_\h - D \Ee)^{-1} C .
\ee
\et
There is a natural generalization of (\ref{eqka3}) to  functions of $d$ variables.
One chooses $d$ Hilbert spaces $\h_1,\dots,\h_d$, lets $\h = \h_1 \oplus \dots \oplus \h_d$,
lets $\Ee = 
\l^1 I_{\h_1} \oplus  \dots \oplus \l^d I_{\h_d} $,
and then, for any isometry $V$ as in (\ref{eqka22}), let
\be
\label{eqka4}
\psi(\lambda) \= A + B \Ee (I_\h - D \Ee)^{-1} C .
\ee
The set of all functions $\psi$ that are realizable
in this way is exactly the Schur-Agler class, 
a class of analytic functions of $d$ variables 
that can also be defined as
\be
\{ \psi \, : \, \| \psi (T_1, \dots, T_d ) \| \leq 1 \quad \forall \ {\rm commuting\ contractive\ matrices\ }
(T_1,\dots,T_d) \}  .
\label{eqka5}
\ee
Von Neumann's inequality \cite{vonN51} is the assertion that for $d=1$ the Schur-Agler class equals 
$\hido$;
And\^o's inequality \cite{and63} is the equality for $d=2$.
Once $d > 2$, the Schur-Agler class is a proper subset of the closed unit ball of $H^\i(\D^d)$ \cite{var74, cradav}.
Many of the results in Sections~\ref{secc} to \ref{secff} are true, with similar proofs, 
for the Schur-Agler class in higher dimensions
(or rather the norm for which this is the unit ball)\footnote{Specifically,
Theorems \ref{thmb3}, \ref{thmke1}, \ref{thmkg2} and $(i) \Leftrightarrow (iii)$ of Theorem \ref{thmhd2}.}, 
but it is not know how to generalize them to
$H^\i(\D^d)$.


The usefulness of the realization formula stems primarily not from its ability to represent 
functions, but to produce functions with desired properties with the aid of a suitably chosen 
isometry $V$ (dubbed a {\em lurking isometry} by Joe Ball).
An example of a lurking isometry argument is the proof of Pick's theorem in Section~\ref{secc}.

It is well-known that equality occurs in the Schwarz lemma on the disk
only for M\"obius transformations.
The Schwarz lemma on $\D^d$ is the following (see \cite{rud69} for a proof). 
\bprop
If $f$ is in $H_1^\i(\D^d)$, then
\be
\label{eqbs}
\sum_{r=1}^d (1 - |\lambda^r|^2 ) \left| \frac{\partial f}{\partial \lambda^r} \, (\lambda) 
\right|
\ \leq \
1 - |f(\lambda)|^2 .
\ee
\eprop
G.~Knese 
\cite{kn07b}
proved that 
equality in (\ref{eqbs}) resulted in a curious form of the realization formula.
(Notice that for $d \geq 3$, the hypothesis is that $f$ lie in $H^\i_1(\D^d)$, but the conclusion
means $f$ must be in the Schur-Agler class.)
\bt
Suppose $f \inn H_1^\i(\D^d)$ is a function of all $d$ variables, and
equality holds in (\ref{eqbs}) everywhere on $\D^d$.
This occurs if and only if $f$ has a representation as in 
(\ref{eqka4})
where each space $\h_r$ is one-dimensional, and the unitary $V$ 
is symmetric (equal to its own transpose).
\et

Analyzing the realization formula, J.M.~Anderson, M.~Dritschel and J.~Rovnyak were able to obtain
the following higher derivative version of the Schwarz lemma \cite{adr08}:
\bt
Let $f$ be in $H_1^\i(\D^2)$, and $n_1, n_2$ non-negative integers with $n = n_1 + n_2$.
Let $\l = (z,w)$ be in $\D^2$, with $|\l| = \max(|z|,|w|)$. Then
$$
\left| \frac{\partial^n f}{\partial^{n_1} z \partial^{n_2} w} \,  
\right|
\ \leq \
(n-2)! 
\,
\frac{1-|f(\l)|^2}{(1-|\l|)^{n-1}} 
\left[
\frac{n_1^2 - n_1}{1 - |z|^2} + 
\frac{2 n_1 n_2 }{\sqrt{1 - |z|^2} \sqrt{1-|w|^2}} + 
\frac{n_2^2 - n_2}{1 - |w|^2}  
\right] .
$$
\et

\section{Pick Problem}
\label{secc}
The Pick problem on the disk is to determine, given $N$ points $\l_1,\dots,\l_N$ in $\D$ and
$N$ complex numbers $w_1, \dots, w_N$, whether there exists $\phi \inn \hido$ such that
$$
\phi(\l_i) = w_i, \qquad i=1,\dots, N.
$$
G. Pick proved \cite{pi16} that the answer is yes if and only if the $N$-by-$N$ matrix
\be
\label{cp}
\left(
\frac{1- w_i \bar w_j}{1-\l_i \bar \l_j}
\right)
\ee
is positive semi-definite.

D.~Sarason realized in \cite{sar67} that Pick's theorem could be proved by showing that
operators that commuted with the backward shift on an invariant subspace could be lifted to
operators that commute with it on all of $H^2$; this result was then generalized by B. Sz-Nagy and C. Foia\c{s}
to the commutant lifting theorem \cite{sznfoi68a}.
Here is a proof of Pick's theorem using a lurking isometry.

\bp
If (\ref{cp}) is of rank $M$, then
one can find vectors $\{ g_i
\}_{i=1}^N$ in $\C^M$ such that
\be
\frac{1- w_i \bar w_j}{1-\l_i \bar \l_j}
\= \la g_i , g_j \ra_{\C^M} .
\label{eqecc2}
\ee
We can rewrite
(\ref{eqecc2}) as
\be
1 + \la \l_i g_i , \l_j g_j \ra \= w_i \bar w_j + 
\la g_i , g_j \ra .
\label{eqecc3}
\ee
The lurking isometry $V : \C \oplus \C^M 
\to \C \oplus \C^M$ is defined by
\be
\label{eqecc4}
V:
\left(
\begin{array}{c}
1 \\
\l_i g_i
\end{array}
\right)
\mapsto
\left(
\begin{array}{c}
w_i\\
g_i
\end{array}
\right).
\ee
We extend linearly to the span of
\be
\label{eqc38}
\left\{
\left(
\begin{array}{c}
1 \\
\l_i g_i
\end{array}
\right)
\ : \ i=1,\dots,N \right\} ,
\ee
and if this is not the whole
space $\C \oplus \C^M$, we extend $V$ arbitrarily so that it remains isometric.
Write $V$ as
$$
V \=
\bordermatrix{&\C &\C^M \cr
\C &A & B \cr
\C^M &C  & D} ,
$$
and define $\phi$ by
\be
\label{eqc39}
\phi(\l) \= A + \l B(I-\l D)^{-1}C.
\ee
By the realization formula Thm.~\ref{thmec1}, $\phi$ 
is in $\hido$.
Moreover, as  (\ref{eqecc4}) implies that 
\beq
A + B \l_i g_i &\=& w_i \\
C + D \l_i g_i &=& g_i ,
\eeq
we get that 
$$
\left( I - \l_i D \right)^{-1} C \= g_i,
$$
and hence
$$
\phi(\l_i) \= A + \l_i B g_i \= w_i,
$$
so $\phi$ interpolates.
\ep
(It is not hard to show that $\phi$ is actually a Blaschke product of degree $M$).

A similar argument using Thm.~\ref{thmka2} solves the Pick problem on the bidisk.
The theorem was first proved, by a different method, 
in \cite{ag1}.

\bt
\label{thmb3}
Given points $\l_1,\dots,\l_N$ in $\Dt$ and complex numbers $w_1,\dots, w_N$, there is a function
$\phi \inn \hibo$ that maps each $\l_i$ to the corresponding $w_i$ if and only if there are 
positive semi-definite matrices $\Gamma^1$ and $\Gamma^2$ such that
\be
\label{eqc6}
1 - w_i \bar w_j \= (1 - \l_i^1 \bar \l_j^1) \Gamma^1_{ij} +
(1 - \l_i^2 \bar \l_j^2) \Gamma^2_{ij} .
\ee
\et

On the polydisk, a necessary condition to solve the Pick problem analagous to (\ref{eqc6}) has recently
been found by A.~Grinshpan, D.~Kaliuzhnyi-Verbovetskyi, V.~Vinnikov and H.~Woerdeman
\cite{gkvw08}. As of this writing, it is unknown if the condition is also sufficient,
but we would conjecture that it is not.
\bt
\label{thmb4}
Given points $\l_1,\dots,\l_N$ in $\D^d$ and complex numbers $w_1,\dots, w_N$, 
a necessary condition for there to be  a function
$\phi \inn \hidd$ that maps each $\l_i$ to the corresponding $w_i$ is:
For every $ 1 \leq p < q \leq d$, there are 
positive semi-definite matrices $\Gamma^p$ and $\Gamma^q$ such that
\be
\label{eqc7}
1 - w_i \bar w_j \= \prod_{r \neq q} (1 - \l_i^r \bar \l_j^r) \Gamma^q_{ij} +
\prod_{r \neq p} (1 - \l_i^r \bar \l_j^r) \Gamma^p_{ij} .
\ee
\et

\section{Nevanlinna Problem}
\label{secd}
If the Pick matrix (\ref{cp})
is singular (\ie if $M < N$)
then the solution is unique; otherwise it is not. R. Nevanlinna found a parametrization
of all solutions in this latter case \cite{nev29} (see also \cite{bh83} for a more modern approach).
\bt
\label{thmd1}
If (\ref{cp}) is invertible, there is a $2$-by-$2$ contractive matrix-valued
function
$$
G \= \left( \begin{array}{cc}
G_{11} & G_{12} \\
G_{21} & G_{22} \end{array}
\right)
$$
such that the set of all solutions of the Pick problem is given by
$$
\{ \phi = G_{11} + G_{12} \frac{\psi G_{21}}{1-G_{22} \psi} \ : \ \psi \inn \hido \} .
$$
\et 

On the bidisk, we shall discuss uniqueness in Section~\ref{sech} below.
Consider now the non-unique case.
Let $\phi$ be in $\hibo$, and so by Theorem~\ref{thmka2} it has a representation
as in (\ref{eqka3}). Define vector-valued functions $F_1, F_2$ by
\be
\label{eqke01}
\bordermatrix{  &\C \cr
\h_1 &F_1(\l)\cr
\h_2 & F_2(\l)}
\ := \
(I_\h - D \Ee)^{-1} C .
\ee
For a given solvable Pick problem with a representation as (\ref{eqc6}), 
say that $\phi$ is affiliated with $(\G^1,\G^2)$ if, for some representation
of $\phi$ and $F_1,F_2$ as in (\ref{eqke01}), 
\beq
F_1(\l_i)^\ast F_1 (\l_j) &\=& \Ga^1_{ij} \\
F_2(\l_i)^\ast F_2 (\l_j) &\=& \Ga^2_{ij}
\eeq
for $i,j = 1, \dots N$.
The situation is complicated by the fact that for a given $\phi$, the pairs
$(\G^1,\G^2)$ with which it is affiliated {\em may or may not} be unique.
J. Ball and T. Trent \cite{baltre98} proved:
\bt
\label{thmke1}
Given a solvable Pick problem, with a representation as in 
(\ref{eqc6}), there is a matrix-valued function $G$
$$
G \= 
\bordermatrix{&\C  &\C^M\cr
\C &G_{11} & G_{12}\cr
\C^M&G_{21} &G_{22}}
$$
in the closed unit ball of $H^\i(\D^2,B(\C \oplus \C^M,
\C \oplus \C^M))$, 
such that
the function $\phi$ solves the Pick problem and is affiliated with
$(\G^1,\G^2)$ 
if and only if it can be written as 
\be
\label{eqke6}
\phi \= G_{11} + G_{12} \Psi (I - G_{22} \Psi)^{-1}
G_{21}
\ee
for some $\Psi$ in $H^\i_1(\D^2,B(\C^M,\C^M))$.
\et

\section{Takagi Problem}
\label{sece}

The case where the Pick matrix (\ref{cp}) has some negative eigenvalues was first studied by
T.~Takagi \cite{tak29}, and later by many other authors
\cite{aak71,nud77,bh83}. See the book \cite{bgr90} for an account.
The principal difference is that if one wishes to interpolate with a unimodular function
(\ie a function that has modulus one on the unit circle $\T$), then one has to allow poles inside $\D$. A typical result
is
\bt
\label{thmaa}
Suppose the Pick matrix is invertible, and has $\pi$ positive eigenvalues and
$\nu$ negative eigenvalues. Then there exists a meromorphic interpolating 
function $\phi$ that is unimodular, and is the quotient of a Blaschke product of degree
$\pi$ by a Blaschke product of degree $\nu$.
\et

If $\G$ is not invertible, the problem is degenerate.
It turns out that there is a big difference between solving the problem of finding
Blaschke products $f,g$ such that
$$
f(\l_i) \= w_i g(\l_i) 
$$
and the problem of solving
$$
f(\l_i) / g(\l_i) \= w_i .$$
(The difference occurs if $f$ and $g$ both vanish at some node $\l_i$; in the first problem
the interpolation condition becomes vacuous, but in the second one needs a relation on the derivatives).
The first problem  is more easily handled as the limit of non-degenerate problems;
see the paper \cite{bkr03} for recent developments on this approach.
The second version of the problem has only recently beeen solved, by H.~Woracek,
using Pontryjagin spaces
\cite{wor08}.

\begin{question}
What is the right version of Theorem~\ref{thmaa} on the bidisk?
\end{question}

\section{Interpolating Sequences}
\label{secf}

Given a sequence $\{ \l_i \}_{i=1}^\i$ in the polydisk $\D^d$, we
say it is interpolating for $H^\i(\D^d)$ if, for any bounded
sequence $\{ w_i \}_{i=1}^\i$, there is a function $\phi$ in 
$H^\i(\D^d)$ satisfying $\phi(\l_i) = w_i$. 
L.~Carleson characterized interpolating sequences on $\D$ in \cite{car58}.

Before stating his theorem, let us introduce some definitions.
A {\em kernel} on $\D^d$ is a positive semi-definite function
$k : \D^d \times \D^d \to \C$, \ie a function such that for any
choice of $\l_1,\dots, \l_N$ in $\D^d$ and any complex numbers
$a_1,\dots, a_N$, we have
$$\sum a_i \bar a_j k(\l_i, \l_j) \ \geq \ 0 .
$$
Given any kernel $k$ 
on $\D^d$, a sequence $\{ \l_i \}_{i=1}^\i$ has an associated
Grammian $G^k$, where
$$
[G^k]_{ij} \= \frac{k(\l_i,\l_j)}{\sqrt{k(\l_i,\l_i) \, k(\l_j,\l_j)}}
.
$$
We think of $G^k$ as an infinite matrix, representing an operator on
$\ell^2$ (that is not necessarily bounded).
When $k$ is the \sz kernel on $\D^d$
\be
\label{eqsz3}
k^S(\z,\l) \= \frac{1}{(1- \z^1 \bar \l^1)(1 - \z^2 \bar \l^2)
\cdots (1- \z^d \bar \l^d)},
\ee
we call the associated Grammian the \emph{\sz Grammian}.
The \sz kernel
is the reproducing kernel for the Hardy space $H^2(\D^d) = P^2(m)$, 
where $m$ is $d$-dimensional Lebesgue
measure on the distinguished boundary $\T^d$ of $\D^d$.

An analogue of the pseudo-hyperbolic metric on the polydisk is the
\emph{Gleason distance},
defined by
$$
\rho(\zeta,\l) \ := \
\sup \{|\phi(\zeta)| : \| \phi \|_{H^\i(\D^d)} \leq 1, \phi(\l) = 0 \} .
$$
We shall call a sequence $\{ \l_i \}_{i=1}^\i$ {\em weakly separated} if there
exists $\vare > 0$ such that, for all $i \neq j$, the Gleason distance
$\rho(\l_i,\l_j) \geq \vare$. 
We call the sequence {\em strongly separated} if there exists $\vare > 0 $ such that,
for all $i$, there is a function $\phi_i$ in $\hido$ such that
$$
\phi_i (\l_j) \=
\left\{ 
\begin{array}{l}
\vare, \qquad j = i \\
0, \qquad j \neq i
\end{array} \right.
$$
In $\D$, a straightforward argument using Blaschke products shows that a sequence is strongly separated
if and only if 
$$
\prod_{j\neq i} \rho(\l_i, \l_j) \geq \vare \qquad \forall\  i .
$$

We can now state Carleson's theorem. Let us note that he proved it using 
function theoretic methods, but later H.~Shapiro and A.~Shields \cite{shashi61}
found a Hilbert
space approach, which has proved to be more easily generalized, \eg to characterizing 
interpolating sequences in the multiplier algebra of the Dirichlet space \cite{marsun}.

\bt
\label{thmf1}
On the unit disk, the following are equivalent:

{\rm (1)} There exists $\vare > 0$ such that
$$
\prod_{j\neq i} \rho(\l_i, \l_j) \geq \vare \qquad \forall\  i .
$$

{\rm (2)} The sequence $\{ \l_i \}_{i=1}^\i$ is an interpolating sequence
for
$H^\i(\D)$.

{\rm (3)} The sequence $\{ \l_i \}_{i=1}^\i$ is weakly separated and the
associated \sz Grammian is
a bounded operator on $\ell^2$.
\et

In 1987 B. Berndtsson, S.-Y. Chang and K.-C. Lin
proved the following theorem \cite{bcl}:
\bt
\label{thmkg1}
Let $d \geq 2$. Consider the three statements

{\rm (1)} There exists $\vare > 0$ such that
$$
\prod_{j\neq i} \rho(\l_i, \l_j) \geq \vare \qquad \forall i.
$$

{\rm (2)} The sequence $\{ \l_i \}_{i=1}^\i$ is an interpolating sequence
for
$H^\i(\D^d)$.

{\rm (3)} The sequence $\{ \l_i \}_{i=1}^\i$ is weakly separated and the
associated \sz Grammian is
a bounded operator on $\ell^2$.

Then {\rm (1)} implies {\rm (2)} and {\rm (2)} implies {\rm (3)}. 
Moreover the converse of
both
these implications is false.
\et

We call the kernel $k$ on $\D^d$ {\em admissible} 
if, for each $1 \leq r \leq d$, the function
\[
( 1 - \z^r \bar \l^r) k(\z,\l) 
\]
is positive semidefinite. (This is the same as saying that multiplication by
each coordinate function on the Hilbert function space with reproducing kernel 
$k$ is a contraction.)

On the unit disk, all admissible kernels are in some sense compressions of
the \sz kernel, and so to prove theorems about $\hid$ one can often just use
the fact that it is the multiplier algebra of $H^2$. 
On the bidisk, there is no single dominant kernel, and one must look at 
a huge family of them.
That is the key idea needed in theorems~\ref{thmka2} and \ref{thmb3}, and it 
allows a different generalization of Theorem~\ref{thmf1}, which was proved in \cite{agmc_isb}.
(If this paragraph seems cryptic, there is a more 
detailed exposition of this point of view
in \cite{ampi}).

For the following theorem, let
$\{e_i\}_{i=1}^\i$ be an orthonormal basis for $\ell^2$.
\bt
\label{thmkg2}
Let $\{ \l_i \}_{i=1}^\i $ be a sequence in $\D^2$. The following
are
equivalent:

{\rm (i)} $\{ \l_i \}_{i=1}^\i $ is an interpolating sequence for $\hb$.

{\rm (ii)} The following two conditions hold.

\noindent
{\rm (a)}  For all admissible kernels $k$,
their normalized Grammians are uniformly bounded:
$$
 G^k  \ \leq \ M I
$$
for some positive constant $M$.

{\rm (b)} For all admissible kernels $k$,
their normalized Grammians are uniformly bounded below:
$$
N G^k \ \geq \  I
$$
for some positive constant $N$.

{\rm (iii)} The sequence  $\{ \l_i \}_{i=1}^\i $ is strongly separated
and
condition {\rm  (a)} alone holds.

{\rm (iv)} Condition {\rm (b)} alone holds.


Moreover, Condition {\rm (a)} is equivalent 
to both {\rm (a${}^\prime$)} and {\rm(a${}^{\prime\prime}$)}:

\noindent
{\rm (a${}^\prime$)}: There exists 
a constant $M$ and positive semi-definite
infinite
matrices $\G^1$ and $\G^2$ such that
$$
M \delta_{ij} - 1 \= \G^1_{ij} (1 - \bar \l_i^1 \l_j^1) + \G^2_{ij}
(1 - \bar \l_i^2 \l_j^2) .
$$

\noindent
{\rm ( a${}^{\prime\prime}$)}: There exists a 
function $\Phi$ in $H^\i(\D^2,B(\ell^2,\C))$ of
norm at most $\sqrt{M}$ such that
$ \Phi(\l_i)  e_i = 1 $.

Condition {\rm (b)} is equivalent to both {\rm (b${}^\prime$)} 
and {\rm (b${}^{\prime\prime}$)}:

\noindent
{\rm (b${}^\prime$)}: There exists a constant $N$ and positive
semi-definite infinite
matrices $\Delta^1$ and $\Delta^2$ such that
$$
N  -   \delta_{ij}  \= \Delta^1_{ij} (1 - \bar \l_i^1 \l_j^1) +
\Delta^2_{ij}
(1 - \bar \l_i^2 \l_j^2) .
$$

\noindent
{\rm (b${}^{\prime\prime}$)}: There 
exists a function $\Psi$ in $H^\i(\D^2,B(\C,\ell^2))$ of
norm at most $\sqrt{N}$ such that
$ \Psi(\l_i) = e_i$.
\et

Neither Theorem~\ref{thmkg1} nor \ref{thmkg2} are fully satisfactory. For example, the following
is still an unsolved problem:

\begin{question}
If a sequence on $\D^2$ is strongly separated, is it an interpolating sequence?
\end{question}

\section{Corona Problem}
\label{secff}

The corona problem on a domain $\Omega$ asks whether, whenever one is given $\phi_1,\dots, \phi_N$
in $H^\i(\Omega)$ satisfying
\be
\label{eqff1}
\sum_{i=1}^N |\phi_i(\l)|^2 \ \geq \ \vare > 0,
\ee
there always exist $\psi_1, \dots, \psi_N $ in $H^\i(\Omega)$ satisfying
\be
\label{eqff2}
\sum_{i=1}^N \phi_i \psi_i = 1 .
\ee
If the answer is affirmative, the domain is said to have no corona.

Carleson proved that the disk has no corona in \cite{car62}.
The most striking example of our ignorance about the bidisk is that the answer there
is still unknown.

\begin{question}
Is the corona theorem true for $\D^2$? 
\end{question}
The best result known is due to T.~Trent \cite{tre06}, who proved that
a solution can be found with the $\psi_i$'s in a specific Orlicz space
${\rm exp}(L^{1/3})$, which is contained in $ \cap_{p < \i} H^p(m)$.

There is a version of the corona theorem, the Toeplitz-corona theorem, proved
at various levels of generality by several authors
\cite{arv75}, \cite{sznfoi76}, \cite{sch78}, \cite{ros80}.
We use $k^S$ as in (\ref{eqsz3}) (with $d=1$).

\bt
\label{thmhd1}
Let $\phi_1,\dots,  \phi_N$ be in $\hid$ and $\delta > 0$. Then the
following are
 equivalent:

{\rm (i)} The function
\be
\label{eqff6}
\left[ \sum_{i=1}^N \phi_i(\z) \overline{ \phi_i (\l)}
-\delta \right] k^S(\z,\l) 
\ee
is positive semi-definite on $\D \times \D$.

{\rm (ii)} The multipliers $M_{\phi_i}$ on $H^2$ satisfy  the inequality
\be
\label{eqhd2}
\sum_{i=1}^N  M_{\phi_i} M_{\phi_i}^\ast  \geq \delta
I .
\ee

{\rm (iii)} There exist functions 
$\psi_1,\dots,  \psi_N$ in $\hid$ such that
\[
\sum_{i=1}^N \psi_i \phi_i  = 1
\]
and
\be
\label{eqhd4}
 \sup_{\l \in \D} \left[ \sum_{i=1}^N |\psi_i(\l)|^2 \right] \leq
\frac{1}{\delta}
 .
\ee
\et

The Toeplitz-corona theorem is often considered a weak version of
the corona theorem, because the proof is easier and the hypothesis
({\ref{eqhd2}) is more stringent than (\ref{eqff1}).
It does, however, have a stronger
conclusion: condition (iii) gives the exact best
bound for the norm of the $\psi_i$'s, whereas the 
corona theorem asserts that if (\ref{eqff1}) holds, 
then (iii) holds for {\it some} $\delta >0 $.
(Moreover, in practice, checking the hypothesis (\ref{eqhd2}) is an eigenvalue
problem, and so quite feasible with polynomial data. Checking
(\ref{eqff1}) is a minimization problem over a function on the disk that
one would expect to have many local minima, even if the $\phi_i$'s are polynomials
of fairly low degree).

The Toeplitz-corona theorem does generalize to the bidisk, but again 
it is not enough to check (\ref{eqff6}) for a single kernel 
(or 
(\ref{eqhd2})
on a single Hilbert function space),
but rather one must find a uniform lower bound that works for all admissible kernels.
For details see \cite{baltre98,agmc_bid}.
\bt
\label{thmhd2}
Let $\phi_1,\dots,  \phi_N$ be in $\hib$ and $\delta > 0$. Then the
following are
 equivalent:

{\rm (i)} The function
\be
\label{eqff67}
\left[ \sum_{i=1}^N \phi_i(\z) \overline{ \phi_i (\l)}
-\delta \right] k(\z,\l) 
\ee
is positive semi-definite for all admissible kernels $k$.

{\rm (ii)} For every measure $\mu$ on $\T^2$, 
the multipliers $M_{\phi_i}$ on $P^2(\mu)$ satisfy  the inequality
\be
\label{eqhd22}
\sum_{i=1}^N  M_{\phi_i} M_{\phi_i}^\ast  \geq \delta
I .
\ee

{\rm (iii)} There exist functions 
$\psi_1,\dots,  \psi_N$ in $\hib$ such that
\[
\sum_{i=1}^N \psi_i \phi_i  = 1
\]
and
\be
\label{eqhd43}
 \sup_{\l \in \D^2} \left[ \sum_{i=1}^N |\psi_i(\l)|^2 \right] \leq
\frac{1}{\delta}
 .
\ee
\et

Although Theorem~\ref{thmhd2} seems to depend on the specific properties
of the bidisk, (indeed, using Theorem~\ref{thmka2} one can prove the equivalence
of (i) and (iii) in the Schur-Agler norm on the polydisk), there is a remarkable generalization
by E.~Amar that applies not only to the polydisk, but to any smooth convex domain \cite{am03}.

\bt
\label{thmhd3}
Let $\O$ be a bounded convex domain in $\C^d$
containing the origin, and assume that either $\O$ is $\D^d$ or its boundary is smooth.
Let $X$ be $\T^d$ in the former case, the boundary of $\O$ in the latter.
Let $\phi_1, \dots, \phi_N$ be in $H^\i(\O)$ and $ \delta > 0$. Then the following are equivalent:

{\rm (i)} There exist functions 
$\psi_1,\dots,  \psi_N$ in $\hid$ such that
\[
\sum_{i=1}^N \psi_i \phi_i  = 1
\]
and
\[
 \sup_{\l \in \O} \left[ \sum_{i=1}^N |\psi_i(\l)|^2 \right] \leq
\frac{1}{\delta}
 .
\]

{\rm (ii)} For every measure $\mu$ on $X$, 
the multipliers $M_{\phi_i}$ on $P^2(\mu)$ satisfy  the inequality
\be
\label{eqhd3}
\sum_{i=1}^N  M_{\phi_i} M_{\phi_i}^\ast  \geq \delta
I .
\ee

{\rm (iii)} For every measure $\mu$ on $X$, 
and every $f$ in $P^2(\mu)$, 
there exist functions 
$\psi_1,\dots,  \psi_N$ in $\ptm$ such that
\[
\sum_{i=1}^N \psi_i \phi_i  = f
\]
and
\be
\label{eqhd45}
 \sum_{i=1}^N \| \psi_i \|^2 \leq
\frac{1}{\delta} \| f \|^2
,
\ee
where the norms on both sides of (\ref{eqhd45}) are in $\ptm$.
\et

In \cite{tw08}, T.~Trent and B.~Wick have shown that in Amar's theorem it is sufficient to consider measures
$\mu$ that are absolutely continuous and whose derivatives are bounded away from zero.

\section{Distinguished and Toral Varieties}
\label{sech}

A Pick problem is called {\em extremal} if it is solvable with a function of norm $1$, but not
with anything smaller. In one dimension, this forces the solution to be unique.
(In the notation of Section~\ref{secc}, this corresponds to the Pick matrix (\ref{cp}) being singular,
the vectors in (\ref{eqc38}) spanning $\C^{1+M}$, and the unique solution being (\ref{eqc39}).)
On the bidisk, problems can be extremal in either one or two dimensions.
For example, consider the problems:
$w_1 = 0, w_2 = 1/2, \, \l_1 = (0,0)$
and $\l_2$ either $(1/2,0)$ or $(1/2,1/2)$.
The first problem has the unique solution $z$; the latter problem has a unique solution
on the one-dimensional set $\{ z = w \}$, but is not unique off this set.
Indeed, Theorem~\ref{thmke1} beomes in this case that the general solution is given by
\[
\phi(z,w) \= t z + (1-t)w + t(1-t)(z - w)^2
\frac{\Psi}{1 - [(1-t)z + t w]\Psi} ,
\]
where $\Psi$ is any function in $H^\i_1(\D^2)$ and $t$ is any number in $[0,1]$.

If an extremal Pick problem on $\D^2$ does not have a solution that is unique on all
$\D^2$, then the set on which it is unique must be a variety\footnote{
We use the word variety where algebraic geometers would say algebraic set -- \ie we do not
require that a variety be irreducible.} (the zero set of a polynomial).
But this is not an arbitrary variety --- it has special properties.

Let $\E$ be the exterior of the closed disk, $\C \setminus \overline{\D}$.
Say a variety $V$  in $\C^2$ is {\em toral} if every irreducible component intersects $\T^2$ in an infinite
set, and say it is {\em distinguished} if  
$$
V \ \subset \ \D^2 \cup \T^2 \cup \E^2 .
$$
Distinguished varieties first 
appeared implicitly in the paper \cite{rud69b} by W.~Rudin,
and later in the operator theoretic context of sharpening And\^o's
inequality for matrices \cite{agmc_dv}; 
they turn out to be intimately connected to function
theory on $\D^2$ (see Theorem~\ref{thmj3}, for example).
Toral varieties are related to inner functions \cite{ams06} and to 
symmetry of a variety with respect to the torus \cite{ams08}.
The uniqueness variety was partially described in \cite{agmc_dv, ams06}:

\bt
\label{thmh1}
The uniqueness set for an extremal Pick problem on $\D^2$ is either all of $\D^2$
or a toral variety. In the latter case, it contains a distinguished variety.
\et

It is perhaps the case that the uniqueness set is all of $\D^2$ whenever the
data is in some sense ``generic'' (see \eg \cite{agmc_three}), but how is that made precise?
\begin{question}
When is the uniqueness set all of $\D^2$?
\end{question}

Distinguished varieties have a determinantal representation.
The following theorem was proved in \cite{agmc_dv}, and, more constructively, in \cite{kn08ua}.
\bt
\label{thmh2}
A variety $V$ is a distinguished variety if and only if there is a pure matrix-valued rational inner
function $\Psi$
on the disk such that
$$
V \cap \D^2 \= \{ (z,w) \inn \D^2 : \det( \Psi(z) - wI) = 0 \} .
$$
\et
Another way to picture distinguished varieties is by taking the Cayley transform of both variables;
then they become varieties in $\C^2$ with the property that when one coordinate is real, so is the other.

\section{Extension property}
\label{secj}

If $G$ is a subset of $\D^2$,
we shall say that a function $f$ defined on $G$ is holomorphic if,
for every point $P$ in $G$, there is an open ball $B(P,\vare)$ in $\D^2$
and a holomorphic function on the ball whose restriction to $G$ is $f$.
Given such a holomorphic function $f$, one can ask whether there is a single function
$F$ on $\D^2$ that extends it, and, if so, whether $F$ can be chosen with additional properties.

H.~Cartan proved that if $G$ is a subvariety, then a global extension $F$ always exists
\cite{car51} (indeed he proved this on any pseudo-convex domain, the bidisk being just a special
case). If $f$ is bounded, one can ask whether one can find an extension $F$ with the same 
sup-norm. If $G$ is an analytic retract of $\D^2$, \ie there is an analytic map
$r: \D^2 \to G$ that is the identity on $G$, then $F = f \circ r$ will work.
(All retracts of $\D^2$ are either singletons, embedded disks, or the whole bidisk \cite{rud69}.)  
It turns out that extending without increasing the norm is only possible for retracts
\cite{agmc_vn}.
\bt
\label{thmj1}
Let $G\subseteq \mathbb{D}^2$ and assume that $G$ is relatively
polynomially convex (i.e. $G^\wedge \cap \mathbb{D}^2 =G$
where $G^\wedge$ denotes the polynomially convex hull of $G$).
If every polynomial $f$ on $G$ has an extension to a function $F$ in $\hib$
of the same norm, then $G$ is
a retract.
\et
Let us remark that although the theorem can be proved without using operator theory, it
was discovered by studying operators that had $G$ as a spectral set.

One can also ask if a bounded function $f$ has a bounded extension $F$, but with a perhaps greater
norm. G. Henkin and P. Polyakov proved that this can always be done if $G$ is a subvariety of
the polydisk that exits transversely \cite{henpol84}. 
In the case of a distinguished variety, Knese showed how to bound the size of the
extension even when there are singularities on $\T^2$ 
\cite{kn08ua,kn08ub} (he also gives a construction of the function $C$ below):
\bt
\label{thmj3}
Let $V$ be a distinguished variety. Then there is a rational function $C(z)$
with no poles in $\D$ such that, for every polynomial $f(z,w)$ there is a rational
function $F$ which agrees with $f$ on $V \cap \D^2$ and satisfies the estimate
$$
| F(z,w) | \ \leq \ |C(z)| \sup_{(z,w) \inn V} |f(z) | .
$$
If $V$ has no singularities on $\T^2$, then $C$ can be taken to be a constant.
\et

\section{Conclusion}
Paul Halmos contributed in many ways to the development of operator theory.
The purpose of this article is to show that recasting many known results about
$\hid$ in terms of operator theory has been extremely fruitful in understanding
$\hib$. So far, however, it has not helped very much in understanding
$H^\i(\B_2)$, where $\B_2$ is the ball in $\C^2$. There is another kernel on the ball,
$$
k(\z,\l) \= \frac{1}{1 - \z^1 \bar \l^1 - \z^2 \bar \l^2},
$$
introduced by S.~Drury \cite{dru78}, and operator theory has been very effective in studying
this kernel \cite{agmc_cnp,agmc_loc,arv98,csw08,grs02,po91}.
\begin{question}
What is the correct Pick theorem on $H^\i(\B_2)$?
\end{question}

\bibliography{references}

\end{document}

%% file: macros_withbirk.tex
\usepackage{latexsym}
\usepackage{amssymb}
\usepackage{euscript}

\def\mcc{M\raise.5ex\hbox{c}C}
\def\mccarthy{M\raise.5ex\hbox{c}Carthy}
\def\Hu{\H}
\def\sz{Szeg\Hu{o} }

\def\eg{{\it e.g. }}
\def\ie{{\it i.e. }}

\def\ptm{P^2(\mu)}

\def\h{{\cal H}}

\def\hid{H^\infty(\D)}
\def\hb{H^\infty(\D^2)}

\def\Ga{\Gamma}

\def\l{\lambda}
\def\z{\zeta}

\def\vare{\varepsilon}

\let\i=\infty

\def\la{\langle}
\def\ra{\rangle}
\def\={\ = \ }

\def\E{E_\l}

\def\Ee{{\cal E}_{\l}}

\def\C{\mathbb C}

\def\T{\mathbb T}
\def\D{\mathbb D}
\def\Dt{{\mathbb D}^2}
\def\B{\mathbb B}

\def\inn{\ \in \ }

\def\be{\setcounter{equation}{\value{theorem}} \begin{equation}}
\def\ee{\end{equation} \addtocounter{theorem}{1}}
\def\beq{\begin{eqnarray*}}
\def\eeq{\end{eqnarray*}}

\def\bp{{\sc Proof: }}
\def\ep{{}{\hfill $\Box$} \vskip 5pt \par}

\def\bl{\begin{lemma}}
\def\el{\end{lemma}}
\def\bt{\begin{theorem}}
\def\et{\end{theorem}}
\def\bprop{\begin{prop}}
\def\eprop{\end{prop}}
\def\bd{\begin{definition}}
\def\ed{\end{definition}}
\def\bexer{\begin{exercise}}
\def\eexer{\end{exercise}}
\def\bfig{\begin{figure}}
\def\efig{\end{figure}}

\newtheorem{theorem}{Theorem}[section]
\newtheorem{prop}[theorem]{Proposition}
\newtheorem{lemma}[theorem]{Lemma}

\newtheorem{question}{Question}
\newtheorem{definition}[theorem]{Definition}

\renewcommand{\theequation}{\thesection.\arabic{equation}}

%% file: Hilbertandbidisk.bbl
\begin{thebibliography}{10}

\bibitem{aak71}
V.M. Adamian, D.Z. Arov, and M.G. Kre{\u{\ii}}n.
\newblock {Analytic properties of Schmidt pairs for a Hankel operator and the
  generalized Schur-Takagi problem}.
\newblock {\em Math. USSR. Sb.}, 15:31--73, 1971.

\bibitem{ag1}
J.~Agler.
\newblock Some interpolation theorems of {Nevanlinna-Pick} type.
\newblock Preprint, 1988.

\bibitem{ag90}
J.~Agler.
\newblock On the representation of certain holomorphic functions defined on a
  polydisc.
\newblock In {\em Operator Theory: Advances and Applications, {Vol. 48}}, pages
  47--66. {Birkh\"auser}, Basel, 1990.

\bibitem{agmc_bid}
J.~Agler and J.E. M\raise.45ex\hbox{c}Carthy.
\newblock {Nevanlinna-Pick} interpolation on the bidisk.
\newblock {\em J. Reine Angew. {M}ath.}, 506:191--204, 1999.

\bibitem{agmc_cnp}
J.~Agler and J.E. M\raise.45ex\hbox{c}Carthy.
\newblock Complete {Nevanlinna-Pick} kernels.
\newblock {\em J. Funct. Anal.}, 175(1):111--124, 2000.

\bibitem{agmc_loc}
J.~Agler and J.E. M\raise.45ex\hbox{c}Carthy.
\newblock {Nevanlinna-Pick} kernels and localization.
\newblock In A.~Gheondea, R.N. Gologan, and D.~Timotin, editors, {\em
  Proceedings of 17th International Conference on Operator Theory at Timisoara,
  1998}, pages 1--20. Theta Foundation, Bucharest, 2000.

\bibitem{agmc_three}
J.~Agler and J.E. M\raise.45ex\hbox{c}Carthy.
\newblock {The three point Pick} problem on the bidisk.
\newblock {\em New York Journal of Mathematics}, 6:227--236, 2000.

\bibitem{agmc_isb}
J.~Agler and J.E. M\raise.45ex\hbox{c}Carthy.
\newblock Interpolating sequences on the bidisk.
\newblock {\em International J. Math.}, 12(9):1103--1114, 2001.

\bibitem{ampi}
J.~Agler and J.E. M\raise.45ex\hbox{c}Carthy.
\newblock {\em Pick Interpolation and Hilbert Function Spaces}.
\newblock American Mathematical Society, Providence, 2002.

\bibitem{agmc_vn}
J.~Agler and J.E. M\raise.45ex\hbox{c}Carthy.
\newblock Norm preserving extensions of holomorphic functions from subvarieties
  of the bidisk.
\newblock {\em Ann. of Math.}, 157(1):289--312, 2003.

\bibitem{agmc_dv}
J.~Agler and J.E. M\raise.45ex\hbox{c}Carthy.
\newblock Distinguished varieties.
\newblock {\em Acta Math.}, 194:133--153, 2005.

\bibitem{ams06}
J.~Agler, J.E. M\raise.45ex\hbox{c}Carthy, and M.~Stankus.
\newblock Toral algebraic sets and function theory on polydisks.
\newblock {\em J. Geom. Anal.}, 16(4):551--562, 2006.

\bibitem{ams08}
J.~Agler, J.E. M\raise.45ex\hbox{c}Carthy, and M.~Stankus.
\newblock Geometry near the torus of zero-sets of holomorphic functions.
\newblock {\em New York J. Math.}, 14:517--538, 2008.

\bibitem{am03}
E.~Amar.
\newblock On the {Toeplitz}-corona problem.
\newblock {\em Publ. Mat.}, 47(2):489--496, 2003.

\bibitem{adr08}
J.M. Anderson, M.~Dritschel, and J.~Rovnyak.
\newblock {Shwarz-Pick} inequalities for the {Schur-Agler} class on the
  polydisk and unit ball.
\newblock {\em Comput. Methods Funct. Theory}, 8:339--361, 2008.

\bibitem{and63}
T.~{And\^o}.
\newblock On a pair of commutative contractions.
\newblock {\em Acta Sci. Math. (Szeged)}, 24:88--90, 1963.

\bibitem{arv75}
W.B. Arveson.
\newblock Interpolation problems in nest algebras.
\newblock {\em J. Funct. Anal.}, 20:208--233, 1975.

\bibitem{arv98}
W.B. Arveson.
\newblock Subalgebras of {C*-algebras III}: {Multivariable} operator theory.
\newblock {\em Acta Math.}, 181:159--228, 1998.

\bibitem{bgr90}
J.A. Ball, I.~Gohberg, and L.~Rodman.
\newblock {\em Interpolation of rational matrix functions}.
\newblock {Birkh\"auser}, Basel, 1990.

\bibitem{bh83}
J.A. Ball and J.W. Helton.
\newblock {A Beurling-Lax theorem for the Lie group $U(m,n)$ which contains
  most classical interpolation theory}.
\newblock {\em Integral Equations and Operator Theory}, 9:107--142, 1983.

\bibitem{baltre98}
J.A. Ball and T.T. Trent.
\newblock Unitary colligations, reproducing kernel {Hilbert} spaces, and
  {Nevanlinna-Pick} interpolation in several variables.
\newblock {\em J. Funct. Anal.}, 197:1--61, 1998.

\bibitem{bcl}
B.~Berndtsson, S.-Y. Chang, and K.-C. Lin.
\newblock Interpolating sequences in the polydisk.
\newblock {\em Trans. Amer. Math. Soc.}, 302:161--169, 1987.

\bibitem{bkr03}
V.~Bolotnikov, A.~Kheifets, and L.~Rodman.
\newblock {Nevanlinna-Pick} interpolation: {Pick} matrices have bounded number
  of negative eigenvalues.
\newblock {\em Proc. Amer. Math. Soc.}, 132:769--780, 2003.

\bibitem{car58}
L.~Carleson.
\newblock An interpolation problem for bounded analytic functions.
\newblock {\em Amer. J. Math.}, 80:921--930, 1958.

\bibitem{car62}
L.~Carleson.
\newblock Interpolations by bounded analytic functions and the corona problem.
\newblock {\em Ann. of Math.}, 76:547--559, 1962.

\bibitem{car51}
H.~Cartan.
\newblock {\em S\'eminaire {Henri Cartan} 1951/2}.
\newblock W.A. Benjamin, New York, 1967.

\bibitem{Co91}
J.B. Conway.
\newblock {\em The Theory of Subnormal Operators}.
\newblock American Mathematical Society, Providence, 1991.

\bibitem{csw08}
S.~Costea, E.T. Sawyer, and B.D. Wick.
\newblock The corona theorem for the {Drury-Arveson Hardy} space and other
  holomorphic {Besov-Sobolev} spaces on the unit ball in {$\C^n$}.
\newblock http://front.math.ucdavis.edu/0811.0627.

\bibitem{cm95}
C.C. Cowen and B.D. MacCluer.
\newblock {\em Composition operators on spaces of analytic functions}.
\newblock CRC Press, Boca Raton, 1995.

\bibitem{cradav}
M.J. Crabb and A.M. Davie.
\newblock {Von Neumann's inequality for Hilbert space operators}.
\newblock {\em Bull. London Math. Soc.}, 7:49--50, 1975.

\bibitem{dru78}
S.W. Drury.
\newblock A generalization of {von Neumann's} inequality to the complex ball.
\newblock {\em Proc. Amer. Math. Soc.}, 68:300--304, 1978.

\bibitem{grs02}
D.~Greene, S.~Richter, and C.~Sundberg.
\newblock The structure of inner multipliers on spaces with complete
  {Nevanlinna Pick} kernels.
\newblock {\em J. Funct. Anal.}, 194:311--331, 2002.

\bibitem{gkvw08}
A.~Grinshpan, D.~Kaliuzhnyi-Verbovetskyi, V.~Vinnikov, and H.~Woerdeman.
\newblock Classes of tuples of commuting contractions satisfying the
  multivariable {von Neumann} inequality.
\newblock {\em J. Funct. Anal.}, 2008.
\newblock in press.

\bibitem{hal}
P.~R. Halmos.
\newblock Normal dilations and extensions of operators.
\newblock {\em Summa Brasil.\ Math.}, 2:125--134, 1950.

\bibitem{henpol84}
G.M. Henkin and P.L. Polyakov.
\newblock Prolongement des fonctions holomorphes bornées d'une sous-variété du
  polydisque.
\newblock {\em Comptes Rendus Acad. Sci. Paris S\'er. I Math.},
  298(10):221--224, 1984.

\bibitem{kn08ua}
G.~Knese.
\newblock Polynomials defining distinguished varieties.
\newblock To appear.

\bibitem{kn08ub}
G.~Knese.
\newblock Polynomials with no zeros on the bidisk.
\newblock To appear.

\bibitem{kn07b}
G.~Knese.
\newblock A {S}chwarz lemma on the polydisk.
\newblock {\em Proc. Amer. Math. Soc.}, 135:2759--2768, 2007.

\bibitem{marsun}
D.~Marshall and C.~Sundberg.
\newblock Interpolating sequences for the multipliers of the {Dirichlet} space.
\newblock Preprint; see
  http://www.math.washington.edu/$\sim$marshall/preprints/preprints.html, 1994.

\bibitem{nev29}
R.~Nevanlinna.
\newblock {\"Uber beschr\"ankte Funktionen}.
\newblock {\em Ann. Acad. Sci. Fenn. Ser. A}, 32(7):7--75, 1929.

\bibitem{nik02}
N.~K. Nikol'ski{\u{\ii}}.
\newblock {\em Operators, functions and systems: An easy reading}.
\newblock AMS, Providence, 2002.

\bibitem{nud77}
A.A. Nudelman.
\newblock On a new type of moment problem.
\newblock {\em Dokl. Akad. Nauk. SSSR.}, 233:5:792--795, 1977.

\bibitem{pel02}
V.V. Peller.
\newblock {\em Hankel operators and their applications}.
\newblock Springer, New York, 2002.

\bibitem{pi16}
G.~Pick.
\newblock {\"U}ber die {Beschr\"ankungen} analytischer {F}unktionen, welche
  durch vorgegebene {F}unktionswerte bewirkt werden.
\newblock {\em Math. Ann.}, 77:7--23, 1916.

\bibitem{po91}
G.~Popescu.
\newblock {Von Neumann inequality for $(B({\cal H})\sp n)\sb 1$}.
\newblock {\em Math. Scand.}, 68:292--304, 1991.

\bibitem{ros80}
M.~Rosenblum.
\newblock A corona theorem for countably many functions.
\newblock {\em Integral Equations and Operator Theory}, 3(1):125--137, 1980.

\bibitem{rud69}
W.~Rudin.
\newblock {\em Function Theory in {Polydiscs}}.
\newblock Benjamin, New York, 1969.

\bibitem{rud69b}
W.~Rudin.
\newblock Pairs of inner functions on finite {Riemann} surfaces.
\newblock {\em Trans. Amer. Math. Soc.}, 140:423--434, 1969.

\bibitem{sar67}
D.~Sarason.
\newblock Generalized interpolation in {$H^\infty$}.
\newblock {\em Trans.\ Amer.\ Math.\ Soc.}, 127:179--203, 1967.

\bibitem{sch78}
C.F. Schubert.
\newblock The corona theorem as an operator theorem.
\newblock {\em Proc. Amer. Math. Soc.}, 69:73--76, 1978.

\bibitem{shashi61}
H.S. Shapiro and A.L. Shields.
\newblock On some interpolation problems for analytic functions.
\newblock {\em Amer. J. Math.}, 83:513--532, 1961.

\bibitem{sznfoi68a}
B.~Szokefalvi-Nagy and C.~Foia\c{s}.
\newblock Commutants de certains op\'erateurs.
\newblock {\em Acta Sci. Math. (Szeged)}, 29:1--17, 1968.

\bibitem{sznfoi76}
B.~Szokefalvi-Nagy and C.~Foia{\c s}.
\newblock {On contractions similar to isometries and Toeplitz operators}.
\newblock {\em Ann. Acad. Sci. Fenn. Ser. AI Math.}, 2:553--564, 1976.

\bibitem{tak29}
T.~Takagi.
\newblock On an algebraic problem related to an analytic theorem of
  {Carath\'eodory and Fejer}.
\newblock {\em Japan J. Math.}, 1:83--93, 1929.

\bibitem{tre06}
T.T. Trent.
\newblock A vector-valued {$H^p$} corona theorem on the polydisk.
\newblock {\em Integral Equations and Operator Theory}, 56:129--149, 2006.

\bibitem{tw08}
T.T. Trent and B.D. Wick.
\newblock Toeplitz corona theorems for the polydisk and the unit ball.
\newblock http://front.math.ucdavis.edu/0806.3428.

\bibitem{var74}
N.Th. Varopoulos.
\newblock On an inequality of von {Neumann} and an application of the metric
  theory of tensor products to operators theory.
\newblock {\em J. Funct. Anal.}, 16:83--100, 1974.

\bibitem{vonN51}
J.~von Neumann.
\newblock {Eine Spektraltheorie f\"ur allgemeine Operatoren eines unit\"aren
  Raumes}.
\newblock {\em Math. Nachr.}, 4:258--281, 1951.

\bibitem{wor08}
H.~Woracek.
\newblock An operator theoretic approach to degenerated {Nevanlinna-Pick}
  interpolation.
\newblock Preprint.

\end{thebibliography}
